# INTEGRAL ESTIMATES FOR SOLUTIONS OF ORDINARY DIFFERENTIAL EQUATIONS WITH A DISCONTINUITY ON A DOMAIN BOUNDARY

N. G. Dokuchaev

We continue to study the functionals of solutions of ordinary differential equations using probability methods. An investigation started in [1]. In this paper we consider the essentially more complicated case in which a functional includes the moment of time when a path first leaves the domain (this moment depends extremely irregularly on the initial data).

## 1. INTRODUCTION

Let $D \subset \mathbf{R}^n$ be a bounded simply connected domain with a $C^1$-smooth boundary $\partial D$ and $T > 0$ be a number. Let us denote $Q = D \times (0,T)$, $\bar{D} = D \cup \partial D$, $\bar{Q} = Q \cup \partial Q$. Consider the equation

$$\frac{dy}{dt}(t) = f(y(t), t). \tag{1.1}$$

Here $f(x,t) : \mathbf{R}^n \times \mathbf{R} \to \mathbf{R}^n$ is a bounded continuous function having continuous first-order partial derivatives with respect to $x$.

We will use $y^{a,s}(t)$ to denote the solution to (1.1) with $t \geq s$ and the initial condition

$$y(s) = a. \tag{1.2}$$

We specify the quantities

$$\tilde{\tau}^{a,s} = \inf\{t : t \geq s,\ y^{a,s}(t) \notin D\}, \quad \tau^{a,s} = \min\{T, \tilde{\tau}^{a,s}\} \tag{1.3}$$

(moments when the trajectory first leaves the domain) for $(x,s) \in Q$.

We will study the quantities

$$U(x,s) = \zeta(y^{x,s}(T))\,\text{Ind}\,\{\tilde{\tau}^{x,s} > T\}, \tag{1.4}$$

$$V(x,s) = \int_s^{\tau^{x,s}} \varphi(y^{x,s}(t), t)\,dt \tag{1.5}$$

for $(x,s) \in Q$. Here $\text{Ind}\,\{\tilde{\tau}^{x,s} > T\} = 1$ if $\tilde{\tau}^{x,s} > T$ and $\text{Ind}\,\{\tilde{\tau}^{x,s} > T\} = 0$ if $T \geq \tilde{\tau}^{x,s}$. The functions $\zeta(\cdot) : \bar{D} \to \mathbf{R}$ and $\varphi(\cdot,\cdot) : \bar{Q} \to \mathbf{R}$ are, generally speaking, discontinuous, unbounded, and even not Borel measurable (they are only required to be Lebesgue measurable and summable to some power).

These properties of the functions $\zeta(\cdot)$ and $\varphi(\cdot)$ make studying the functionals in (1.4) and (1.5) difficult even for the simple case of $D = \mathbf{R}^n$, $\tau^{x,s} \equiv T$, which was considered in [1]. When $D \neq \mathbf{R}^n$, the problem is even more complicated: even for smooth $f$, $\zeta$, $\varphi$, and $\partial D$ the functionals in (1.4) and (1.5) do not depend on $(x,t)$ continuously. (An example can easily be constructed: let $n = 1$, $D = (-1,1)$ be an interval, $f(x,t) = \cos t$, $T \geq 2\pi$, then $y^{x,s}(t) = x + \sin t - \sin s$, $\tau^{0,0} = \pi/2$, $\tau^{x,0} \to 3\pi/2$ as $x \to -0$). At the same time, quantities like those defined in (1.4) and (1.5) are important numerical characteristics of dynamic processes and are often considered in applications. Therefore, it is interesting to clear up how regularly they depend on $x$, $s$, $\varphi$, and $\zeta$. In this paper, we derive estimates for the $L_2$-norms of the functions $U(\cdot,s)$ and $V(\cdot,s)$ and study their continuity in $L_2(D)$ with respect to $s$.

We will use probability methods; i.e., instead of (1.1) we consider the Itô equations

$$dy_\varepsilon(t) = f(y_\varepsilon(t), t)\,dt + \varepsilon\,dw(t), \tag{1.6}$$

 

which are close to it (for small $\varepsilon$) in some sense. Here $\varepsilon > 0$ is a number and $w(t)$ is the standard $n$-dimensional Wiener process, $t \geq 0$. We will study the quantities

$$U_\varepsilon(x,s) = M\{\zeta(y_\varepsilon^{x,s}(T))\text{Ind}\{\tilde{\tau}_\varepsilon^{x,s} > T\}\}, \tag{1.7}$$

$$V_\varepsilon(x,s) = M \int_s^{\tau_\varepsilon^{x,s}} \varphi(y_\varepsilon^{x,s}(t), t) dt \tag{1.8}$$

for $(x,s) \in \bar{Q}$ and a random solution $y_\varepsilon(t) = y_\varepsilon^{x,s}(t)$ of (1.6) with the initial condition $y(s) = x$, instead of the quantities defined in (1.4) and (1.5). Here $\tilde{\tau}_\varepsilon^{x,s}$ and $\tau_\varepsilon^{x,s}$ are random and are the moments the trajectory leaves the domain for the process $y^{x,s}(t) = y_\varepsilon^{x,s}(t)$ according to (1.3), where $M$ denotes mathematical expectation.

If $\varepsilon > 0$ and $q > 1$, we have $U_\varepsilon \in W_q^{2,1}(Q)$, $V_\varepsilon \in W_q^{2,1}(Q)$ for $\zeta \in L_q(D)$, $\varphi \in L_q(Q)$. This follows from the standard representation [2] of functions in (1.7) and (1.8) as solutions of the first boundary-value problems in $Q$ and from the theorems in [3]. Hence, means shown in (1.7) and (1.8) depend continuously on $(x,s)$ (and even more regularly). Moreover, when $\varepsilon > 0$ sample values of the random quantities $\tau_\varepsilon^{x,s}$ as well as the mean, depend on $(x,s)$ regularly. Later we shall take the limit as $\varepsilon \to 0$ and see that certain properties of $U_\varepsilon$ and $V_\varepsilon$ remain for $U$ and $V$.

In addition to the quantities (1.4) and (1.5), we study the evolution of distribution densities for $y(t)$. Let $a$ in (1.2) be a random vector such that $a \in D$ almost certainly.

Let $a$ have a distribution density $\varrho_s(x) \in L_2(D)$. A process $y^{a,s}(t)$ which breaks when leaving $D$ and satisfies (1.1) and (1.2) is shown to have a density $p(x,t) \in L_2(D)$ ($\forall t \in [s,T]$). This means that

$$P(y^{a,s}(t) \in B, \tau^{a,s} > t) = \int_B p(x,t) dx$$

for any Borel set $B \subset D$ (from here on $P$ denotes the probability of a random event).

Let us introduce some notation. The Hilbert space of a function $u : [s,T] \to H^k$ with finite norm $\|u\|_{X^k[s,T]} = \left( \int_s^T \|u(\cdot,t)\|_{H^k}^2 dt \right)^{1/2}$ is denoted $X^k[s,T]$, $k = 0, 1$. Here $H^0 = L_2(D)$, $H^1 = \mathring{W}_2^1(D)$ are closures in the usual norm of the Sobolev space $W_2^1(D)$ of the sets of finite and smooth functions in $D$, respectively. A Banach space of the function $u : [s,T] \to H^k$ has the finite norm

$$\|u\|_{X_1^k[s,T]} = \int_s^T \|u(\cdot,t)\|_{H^k} dt$$

and is denoted $X_1^k[s,T]$.

A linear normed space of a bounded function $u : [s,T] \to L_2(D)$ equipped with a finite norm

$$\|u\|_{\mathscr{B}[s,T]} = \sup_{t \in [s,T]} \|u(\cdot,t)\|_{L_2(D)}$$

is denoted $\mathscr{B}[s,T]$.

Finally, we use $\mathscr{C}^k[s,T] = C([s,T] \to H^k)$.

We can permit some liberties in the presentation and refer to not only equivalence classes but also Lebesgue measurable functions $u : \bar{D} \times [s,T] \to \mathbf{R}$ as elements of these spaces. For brevity, we denote $X^k = X^k[0,T]$, $X_1^k = X_1^k[0,T]$, $\mathscr{B} = \mathscr{B}[0,T]$, and $\mathscr{C}^k = \mathscr{C}^k[0,T]$.

From here on we denote the norm for a normed space $\mathscr{X}$ by $\|\cdot\|_\mathscr{X}$ and the inner product for a Hilbert space $\mathscr{X}$ by $(\cdot,\cdot)_\mathscr{X}$.

Let us introduce the functions

$$\delta_j(t) = \frac{1}{2} \sup_{x \in D} \left| \sum_{i=1}^n \frac{\partial f_i}{\partial x_i}(x,t) \right|, \quad c_f(s,t) = \exp \int_s^t \delta_f(\varrho) d\varrho, \tag{1.9}$$

where $x_i$ and $f_i$ are components of the vectors $x$ and $f$, respectively.

## 2. EVOLUTION OF A DISTRIBUTION DENSITY

For $\varepsilon > 0$ let us consider the boundary value problems

$$\frac{\partial u_\varepsilon}{\partial t}(x,t) + \sum_{i=1}^{n} \frac{\partial u_\varepsilon}{\partial x_i}(x,t) f_i(x,t) + \frac{\varepsilon^2}{2}\Delta u_\varepsilon(x,t) = 0, \quad (2.1)$$

$$u_\varepsilon(x,T) = \zeta(x), \quad u_\varepsilon(x,t)|_{x \in \partial D} = 0 \quad (2.2)$$

and

$$\frac{\partial p_\varepsilon}{\partial t}(x,t) = -\sum_{i=1}^{n} \frac{\partial}{\partial x_i}(f_i(x,t)p_\varepsilon(x,t)) + \frac{\varepsilon^2}{2}\Delta p_\varepsilon(x,t), \quad (2.3)$$

$$p_\varepsilon(x,s) = \varrho_s(x), \quad p_\varepsilon(x,t)|_{x \in \partial D} = 0 \quad (2.4)$$

in the cylinders $\bar{Q}_s = D \times [s,T]$. Here $\Delta$ denotes the $n$-dimensional Laplacian operator, $x_i$ and $f_i$ are components of the vectors $x$ and $f$, respectively.

The problems (2.1), (2.2) and (2.3), (2.4) can be resolved in terms of

$$u_\varepsilon \in \mathscr{C}^0[s,T] \cap X^1[s,T], \quad p_\varepsilon \in \mathscr{C}^0[s,T] \cap X^1[s,T]$$

for any $\varrho_s \in L_2(D)$, $\zeta \in L_2(D)$ [6, p. 103]. Let us introduce the notation $u_\varepsilon(\cdot,t) = \mathscr{L}_{t,T}(\varepsilon)\zeta$, $p_\varepsilon(\cdot,t) = \mathscr{L}^*_{s,T}(\varepsilon)\varrho$, $p_\varepsilon(\cdot,\cdot) = \mathscr{L}^*_s(\varepsilon)\varrho_s$ for $s \leq t$. The linear operators $\mathscr{L}_{t,T}(\varepsilon) : L_2(D) \to L_2(D)$, $\mathscr{L}^*_{s,T}(\varepsilon) : L_2(D) \to L_2(D)$, and $\mathscr{L}^*_s(\varepsilon) : L_2(D) \to \mathscr{C}^0[s,T]$ can be shown [6, 7] to be continuous and their norms to be majorized by the quantities $c_f(t,T)$, $c_f(s,T)$, and $c_f(s,T)$ as defined according to (1.9), respectively (see the proof of Lemma (1.1) in [1]).

The following lemma can easily be proved by using methods as in [7, Sec. IV].

**Lemma 2.1.** *a) The estimate*

$$\|p_\varepsilon\|_{\mathscr{B}[s,T]} \leq c_0 \|\varrho_s\|_{L_2(D)}, \quad (2.5)$$

*where $c_0 = c_f(0,T)$ is specified by (1.9), is satisfied for $\varepsilon > 0$, $\varrho_s \in L_2(D)$, and $p_\varepsilon = \mathscr{L}^*_s(\varepsilon)\varrho_s$;*

*b) for any $t \in [s,T]$ there exists a sequence $\varepsilon_i = \varepsilon_i(t)$, $i = 0, 1, \ldots$, such that $\varepsilon_i(t) \to 0$, $\varepsilon_i(t) > 0$, the sequence $\{p_{\varepsilon_i}(\cdot,t)\}_{i=1}^{\infty}$ weakly converges in $L_2(D)$ to an element $p(\cdot,t) \in L_2(D)$, and*

$$\|p\|_{\mathscr{B}[s,T]} \leq c_0 \|\varrho_s\|_{L_2(D)}. \quad (2.6)$$

Statement a) of the following theorem is well-known, but we present here a probability proof of it [this proof turns out to be connected with statement b) of the theorem and is of special interest].

**Theorem 2.1.** *a) The function $p(x,t)$ in Lemma 2.1 is uniquely defined as an element of $\mathscr{B}[s,T]$;*

*b) if $\varrho_s \in L_2(D)$ in (2.4) is a distribution density of a random vector $a$ such that $a \in D$ almost certainly, then $p(x,t)$ is the distribution density of a solution $y^{a,s}(t)$ which has a discontinuity on the boundary of the domain, i.e.,*

$$P(y^{a,s}(t) \in B, \tau^{a,s} > t) = \int_B p(x,t)dx \quad (2.7)$$

*is satisfied for any Borel set $B \subset D$ (here $P$ denotes the probability of a random event).*

Now let us introduce the continuous operators $\mathscr{L}^*_s : L_2(D) \to \mathscr{B}[s,T]$ and $\mathscr{L}^*_{s,t} : L_2(D) \to L_2(D)$ such that $p = \mathscr{L}^*_s\varrho_s$, $p(\cdot,t) = \mathscr{L}^*_{s,t}\varrho_s$, where $p$ is the function from Lemma 2.1 and $\mathscr{L}_{s,t}$ denotes the conjugate operator $\mathscr{L}^*_{s,t} : L_2(D) \to L_2(D)$.

The function $p(\cdot,t)$ is weakly continuous in $L_2(D)$ with respect to $t \in [s,T]$. We need a stronger version of this statement.

**Lemma 2.2.** We use $S^+ = \{(s,t) \in [0,T] \times [0,T] : s \leq t\}$. For any set $\{u_{s,t}\} \subset L_2(D)$ such that $(s,t) \in S^+$ and
$$\|u_{s,t}\|_{L_2(D)} \leq \text{const} \quad (\forall (s,t) \in S^+),$$
we have
$$(\mathscr{L}_{s,t}u_{s,t} - u_{s,t}, \xi)_{L_2(D)} \to 0, \quad (\xi, \mathscr{L}^*_{s,t}u_{s,t} - u_{s,t})_{L_2(D)} \to 0$$
for any $\xi \in L_2(D)$ as $t - s \to +0$.

**Theorem 2.2.** *The function $p(\cdot, t)$ in Lemma 2.1 is weakly continuous and strongly right-continuous with respect to $t$ in $L_2(D)$, i.e., $\|p(\cdot, t) - p(\cdot, \theta)\|_{L_2(D)} \to 0$ as $t \to \theta + 0$ for $0 \leq s \leq \theta < t \leq T$.*

*Proof of Theorem 2.1.* It is enough to consider the case in which $\varrho_s(x) \geq 0$ and $\|\varrho_s\|_{L_1(D)} = 1$ (it is clear how to generalize the case because (2.3) and (2.4) is linear). Such a function $\varrho_s(x)$ is the distribution density of a random vector $a$ with values in $D$. We know a solution $p_\varepsilon(\cdot, t)$ of the problem (2.3) and (2.4) to be a distribution density of a solution $y_\varepsilon^{a,s}(t)$ of (1.6) with the boundary condition $y_\varepsilon^{a,s}(s) = a$, which has a discontinuity on the boundary $\partial D$, i.e.,
$$P(y_\varepsilon^{a,s}(t) \in B, \tau_\varepsilon^{a,s} > t) = \int_B p_\varepsilon(x,t)dx$$
for any Borel set $B \subset D$. As before we use $\tau_\varepsilon^{a,s}$ to denote the moment (1.3) when the trajectory first leaves the domain for a process $y_\varepsilon^{a,s}(t)$ satisfying (1.6) and (1.2). Besides, we use $y^{a,s}(t)$ to denote a solution of (1.1) and (1.2) and $\tau^{a,s}$ to denote the corresponding moment (1.3) when the trajectory leaves first the domain. For $\bar{D} = D \cup \partial D$ let us set
$$\tilde\tau^{x,s} = \inf\{t : t \geq s, \ y^{x,s}(t) \notin \bar{D}\}. \tag{2.8}$$

**Proposition 2.1.** *The Lebesgue measure of $x \in D$ for which $\tilde\tau^{x,s} \neq \bar\tau^{x,s}$ is zero (or*
$$P(\tilde\tau^{a,s} \neq \bar\tau^{a,s}) = 0, \tag{2.9}$$
*which is the same if $a$ is a random vector with a distribution density).*

The idea behind the following proof (as well as the complete proof for autonomous systems) was suggested by S. Yu. Pilyugin.

*Proof.* Let us use $\mathscr{T}_x$ to denote a tangent hyperplane to the boundary $\partial D$ at a point $x \in \partial D$. Let $S_0 \subset \partial D \times (-T, T)$ be an open connected $n$-dimensional manifold, $(x,t) \in S_0$, and the surface $S_0$ be specified in a neighborhood of $(x,t)$ by the formula $(x', t') = (\sigma(\xi'), t')$, where $\sigma : \mathbf{R}^{n-1} \to \mathbf{R}$ is a continuously differentiable function, $t' \in (s_1, s_2)$, $s_1 < s_2$.

a) For simplicity, let $f(x,t) \equiv f(x)$. Let us consider a mapping $F : S_0 \to \mathbf{R}^n$ of the form $F(x,t) = y^{x,0}(t)$. For $(x,0) \in S_0$ we use $\Phi(t, \varrho)$ to denote a fundamental matrix of the equation
$$\frac{d\eta(t)}{dt} = \frac{\partial f}{\partial x}(y^{x,0}(t))\eta(t).$$

Clearly, the functions
$$\eta_1(t) = \frac{\partial y^{x,0}(t)}{\partial x}, \quad \eta_2(t) = \frac{\partial y^{x,0}}{\partial t}(t) = f(y^{x,0}(t))$$
satisfy this equation.

It can easily be seen that the differential
$$\frac{dF}{d(x,t)}(x,t) = \left|\frac{\partial y^{x,0}}{\partial x}(t)\frac{d\sigma}{d\xi}(\xi_0), \ f(y^{x,0}(t))\right| = \left|\Phi(t,0)\frac{d\sigma}{d\xi}(\xi_0), \ \Phi(t,0)f(x)\right| \tag{2.10}$$

exists at a point $(x, t) = (\sigma(\xi_0, t))$, $\xi_0 \in \mathbf{R}^{n-1}$. If $f(x) \in \mathscr{T}_x$, then the rank of the matrix in (2.10) is less than $n$. The mapping $F$ is $C^1$-smooth, hence, in view of the Sard theorem [8, p. 58],

$$\text{mes } \{\cup y^{x,0}(t) : (x, t) \in S_0, \ f(x) \in \mathscr{T}_x\} = 0.$$

It can easily be seen that when $f(x, t) \equiv f(x)$, Proposition 2.1 follows from this because if $\bar{\tau}^{x,s} \neq \tilde{\tau}^{x,s}$, then $f(\bar{x}, \tilde{\tau}^{x,s}) \in \mathscr{T}_{\bar{x}}$ for $\bar{x} = y^{x,s}(\tilde{\tau}^{x,s})$.

b) Let us consider the general case. As above, we introduce an $n$-dimensional manifold $S_0$. Let us consider a mapping $\hat{F} : S_0 \times (-T, T) \to \mathbf{R}^{n+1}$ of the form $\hat{F}(z, t) = Y^{z,s}$, where $z = (x, s) \in S_0$, $Y^{z,s}(t)$ is a solution of equations

$$\frac{dY^{z,s}}{dt}(t) = \hat{f}(Y^{z,s}(t)), \quad Y^{z,s}(s) = z$$

(here $\hat{f}(Y) = \text{col}[f(Y_1, \ldots, Y_n, Y_{n+1}), 1]$ for $Y = \text{col}(Y_1, \ldots, Y_n, Y_{n+1})$).

For $x = \sigma(\xi_0)$, $z = (x, s)$ we have $F(z, t) = \text{col}[y^{x,s}(t), t]$,

$$\frac{d\hat{F}}{d(z, t)}(z, t) = \begin{bmatrix} \frac{\partial y^{x,s}}{\partial x}(t)\frac{d\sigma}{d\xi}(\xi_0) & \frac{\partial y^{x,s}}{\partial s}(t) & f(y^{x,s}(t), t) \\ 0_{1,n-1} & 0_{1,1} & 1 \end{bmatrix} \quad (2.11)$$

(here $0_{k,m}$ is a zero in $\mathbf{R}^{k,m}$).

Let $\hat{\Phi}(t, \varrho)$ be the fundamental matrix of the equation

$$\frac{\partial \hat{\eta}}{\partial t}(t) = \frac{\partial f}{\partial x}(Y^{z,s}(t))\hat{\eta}(t),$$

which is specified for $z = (x, s)$. It can easily be seen that

$$\frac{\partial y^{x,s}}{\partial x}(t) = \hat{\Phi}(t, s), \quad \frac{\partial y^{x,s}}{\partial s}(t) = -\hat{\Phi}(t, s)f(z).$$

Substituting these into (2.11) and premultiplying the columns on the right-hand side by $\hat{\Phi}(t, s)^{-1}$, we get

$$\text{rank } \frac{d\hat{F}}{d(\hat{z}, t)}(\hat{z}, t) = 0 \quad \text{for } f(z) \in \mathscr{T}_x.$$

The mapping $\hat{F}$ is $C^1$-smooth, and hence, from the Sard theorem we have the equality for the $(n+1)$-dimensional measure

$$\text{mes } \{\cup \text{col } [y^{x,s}(t), t] : z = (x, s) \in \hat{S}_0, \ f(x, s) \in \mathscr{T}_x, \ t \in (-T, T)\} = 0.$$

This yields Proposition 2.1.

**Proposition 2.2.** *There exists a sequence $\{\varepsilon_i\}_{i=1}^{+\infty}$ such that $\varepsilon_i > 0$, $\varepsilon_i \to +0$ as $i \to +\infty$, and for any domain $B \subset D$ and any $t \in [s, T]$*

$$P(y_{\varepsilon_i}^{a,s}(t) \in B, \ \tau_{\varepsilon_i}^{a,s} > t) \to P(y^{a,s}(t) \in B, \ \tau^{a,s} > t)$$

*as $i \to +\infty$.*

*Proof.* We set $z_\varepsilon(t) = \sup_{\theta \in [s,T]} |y_\varepsilon^{a,s}(\theta) - y^{a,s}(\theta)|^2$ and get

$$Mz_\varepsilon(t) \leq c_1 \{M \sup_{s \leq \theta \leq t} \int_s^\theta |f(y_\varepsilon^{a,s}(t_1), t_1) - f(y^{a,s}(t_1), t_1)|^2 dt_1$$

$$+ M \sup_{s \leq \theta \leq t} |\varepsilon \int_s^\theta dw(t_1)|^2 \} \leq c_2 \int_s^t Mz_\varepsilon(t_1) dt_1 + c_3 \varepsilon^2 (t-s)^2.$$

Here $c_i > 0$ are constants; the second term in the middle of this chain of inequalities is estimated by means of inequality (3) in [9, p. 18]. The Gronwall lemma yields $Mz_\varepsilon(T) \to 0$ as $\varepsilon \to 0$. Hence, there exists a sequence $\{\varepsilon_i\}_{i=1}^{+\infty}$ such that $\varepsilon_i > 0$, $\varepsilon_i \to +0$ as $i \to +\infty$, and

$$z_{\varepsilon_i}(T) = \sup_{t \in [s,T]} |y_{\varepsilon_i}^{a,s}(t) - y^{a,s}(t)|^2 \to 0 \quad \text{as} \quad i \to +\infty$$

is satisfied with probability 1.

From (2.9), we find that

$$\text{Ind}\{y_{\varepsilon_i}^{a,s}(t) \in B, \tau_{\varepsilon_i}^{a,s} > t\} \to \text{Ind}\{y^{a,s}(t) \in \dot{B}, \tau^{a,s} > t\}$$

with probability 1 as $i \to +\infty$ (here Ind$\Gamma$ denotes the indicator function of a random event $\Gamma$). From this and the Lebesgue theorem on limits in integrals, we obtain Proposition 2.2.

We return to Theorem 2.1. For any domain $B \subset D$ and any $t \in [s,T]$, by virtue of Lemma 2.1, there exists a sequence $\{\nu_k\}_{k=1}^{+\infty}$ such that $\nu_k = \varepsilon_{i_k}(t) \to +0$ as $k \to +\infty$, and

$$P(y_{\nu_k}^{a,s}(t) \in B, \tau_{\nu_k}^{a,s} > t) = \int_B p_{\nu_k}(x,t) dx \to \int_B p(x,t) dx$$

at the limit (a subsequence is such that $p(\cdot, t)$ is a weak limit of $p_{\nu_k}(\cdot, t)$ in $L_2(D)$). This and Proposition 2.2 yield Theorem 2.1.

*Remark 2.1.* It follows from Proposition 2.1 that Theorem 2.1 remains valid if we set $\tau^{a,s} = \min\{\bar{\tau}^{a,s}, T\}$ as the moment when the trajectory first leaves a closed set (compare with (1.3), where $\tau^{x,s} = \min\{\tilde{\tau}^{x,s}, T\}$, and $\tilde{\tau}^{x,s}$ is the moment when the trajectory leaves an open set). Moreover, all the assertions in this paper remain valid if we overdetermine $\tau^{x,s}$ and $\tau^{a,s}$ everywhere in the same manner and replace $\tilde{\tau}^{x,s}$ in (1.4) by $\bar{\tau}^{x,s}$. In doing so, the values of the functions in (1.4) and (1.5) will only change for a set $x$ of zero measure.

*Proof of Lemma 2.2.* Let $p_\varepsilon = \mathscr{L}_s^*(\varepsilon) u_{s,t}$, $\xi(\cdot) : D \to \mathbf{R}$ be a smooth function with its support inside $D$. We have

$$(\xi, \mathscr{L}_{s,t}^*(\varepsilon) u_{s,t} - u_{s,t})_{L_2(D)} = \int_s^t \left\{ (p_\varepsilon(\cdot, r), \sum_{i=1}^n \frac{\partial \xi}{\partial x_i}(\cdot) f_i(\cdot, r)) + \frac{\varepsilon^2}{2} (p_\varepsilon(\cdot, r), \Delta \xi(\cdot))_{L_2(D)} \right\} dr$$

Thus,

$$|(\xi, p_\varepsilon(\cdot, t) - u_{s,t})_{L_2(D)}| \leq (t-s) \|p_\varepsilon\|_{\mathscr{B}[s,T]} \|\nabla \xi f + \frac{\varepsilon^2}{2} \Delta \xi\|_{\mathscr{B}[s,T]} \leq (t-s) \|u_{s,t}\|_{L_2(D)} \|\nabla \xi f + \frac{\varepsilon^2}{2} \Delta \xi\|_{\mathscr{B}[s,T]} \to 0.$$

Since smooth functions $\xi$ with a support in $D$ are dense in $L_2(D)$, we get $(u, \mathscr{L}_{s,t}^* u_{s,t} - u_{s,t})_{L_2(D)} \to 0$ ($\forall u \in L_2(D)$). The second statement is proved in a similar manner.

*Proof of Theorem 2.2.* We have $\mathscr{L}_{s,t}^* = \mathscr{L}_{\theta,t}^* \mathscr{L}_{s,\theta}^*$ for $s \leq \theta < t < T$. Hence, we get

$$(\xi, p(\cdot, t) - p(\cdot, \theta))_{L_2(D)} = (\xi, \mathscr{L}_{s,t}^* \varrho_s - \mathscr{L}_{s,\theta}^* \varrho_s)_{L_2(D)} = (\xi, (\mathscr{L}_{\theta,t}^* - I) \mathscr{L}_{s,\theta}^* \varrho_s)_{L_2(D)} \to 0 \quad \text{as} \quad t - \theta \to 0$$

for smooth functions $\xi$ with support in $D$ by virtue of Lemma 2.2 and because $\|\mathscr{L}_{s,\theta}^* \varrho\|_{L_2(D)}$ is bounded [here $I : L_2(D) \to L_2(D)$ is the identity operator, $Iu \equiv u$]. Weak continuity is proved.

Now let us prove $p(\cdot,t)$ to be strongly continuous on one side in $L_2(D)$. For $s \leq \theta < t \leq T$ denote $r = p(\cdot,t) - p(\cdot,\theta)$. We get

$$\|r\|^2_{L_2(D)} = (\mathscr{L}^*_{s,t}\varrho - \mathscr{L}^*_{s,\theta}\varrho,\ \mathscr{L}^*_{s,t}\varrho - \mathscr{L}^*_{s,\theta}\varrho)_{L_2(D)}$$
$$= (\varrho, (\mathscr{L}_{s,t} - \mathscr{L}_{s,\theta})(\mathscr{L}^*_{s,t} - \mathscr{L}^*_{s,\theta})\varrho)_{L_2(D)}$$
$$= (\varrho, \mathscr{L}_{s,\theta}(\mathscr{L}_{\theta,t} - I)(\mathscr{L}^*_{\theta,t} - I)\mathscr{L}^*_{s,\theta}\varrho)_{L_2(D)} \to 0$$

as $t \to \theta + 0$ by virtue of Lemma 2.2. This completes the proof of Theorem 2.2.

## 3. CONTINUITY OF SOLUTIONS OF CONJUGATE PROBLEMS

Let $\mathscr{L}_s : X^0[s,T] \to L_2(D)$ be a conjugate operator $\mathscr{L}^*_s : L_2(D) \to X^0[s,T]$. Later we denote the restriction of $\varphi : \bar{D} \times [0,T] \to \mathbf{R}$ to $D \times [s,T]$ by $\varphi|_{[s,T]}$ for $\varphi \in X^0_1$.

The definition of the operator $\mathscr{L}_s$ yields

$$(v(\cdot,s), \varrho_s(\cdot,s))_{L_2(D)} = \int_s^T (\varphi(\cdot,t), p(\cdot,t))_{L_2(D)} dt \qquad (3.1)$$

for $\varphi \in X^0$, $p = \mathscr{L}^*_s\varrho_s$, $\varrho_s \in L_2(D)$ and the function $v(x,s) = \mathscr{L}_s(\varphi|_{[s,T]})$. We use the notation $v = L\varphi$ for the functions $v(x,s) = \mathscr{L}_s(\varphi|_{[s,T]})$.

**Lemma 3.1.** *The operator $L$ maps $X^0$ into the space $\mathscr{B}$ and can be extended to be a continuous operator $L : X^0_1 \to \mathscr{B}$. The norm of the operator $L : X^0_1 \to \mathscr{B}$ is lower than or equal to $c_0 = c_f(0,T)$, which is specified by (1.9).*

**Theorem 3.1.** *Let $\zeta \in L_2(D)$, $\varphi \in X^0_1$, $u(x,t) = \mathscr{L}_{t,T}\zeta$, $v(x,t) = L\varphi$. Then $u \in \mathscr{B}$, $v \in \mathscr{B}$, the functions $u(\cdot,t)$ and $v(\cdot,t)$ are weakly continuous with respect to $t$ in $L_2(D)$ and strongly left-continuous in $L_2(D)$ (i.e., $\|u(\cdot,t) - u(\cdot,\theta)\|_{L_2(D)} \to 0$, $\|v(\cdot,t) - v(\cdot,\theta)\|_{L_2(D)} \to 0$ as $\theta \to t - 0$).*

*Proof of Lemma 3.1.* Formula (3.1) shows that

$$\|v(\cdot,s)\|_{L_2(D)} = \sup_{\varrho \in \mathfrak{B}}(v(\cdot,s),\varrho)_{L_2(D)} = \sup_{\varrho \in \mathfrak{B}}(\varphi, \mathscr{L}^*_s\varrho)_{X^0[s,t]} \leq \|\varphi\|_{X^0_1[s,t]} \sup_{\varrho \in \mathfrak{B}}\|\mathscr{L}^*_s\varrho\|_{\mathscr{B}[s,t]} \leq c_0\|\varphi\|_{X^0_1},$$

where $\mathfrak{B} = \{\varrho \in L_2(D) : \|\varrho\|_{L_2(D)} \leq 1\}$, $c_0 = c_f(0,T)$ is the constant introduced in Sec. 2 and is common for all $s \in [0,T]$, $\varphi \in X^0$. Applying this estimate to $\varphi \in X^0_1$, we get Lemma 3.1.

*Proof of Theorem 3.1.* Let $s \leq \theta < t \leq T$. For $\xi \in L_2(D)$ we have

$$(\xi, u(\cdot,t) - u(\cdot,\theta))_{L_2(D)} = (\xi, \mathscr{L}_{t,T}\zeta - \mathscr{L}_{\theta,T}\zeta)_{L_2(D)}$$
$$= (\xi, \mathscr{L}_{t,T}\zeta - \mathscr{L}_{\theta,t}\mathscr{L}_{t,T}\zeta)_{L_2(D)} = (\xi, (I - \mathscr{L}_{\theta,t})\mathscr{L}_{t,T}\zeta)_{L_2(D)} \to 0$$

as $t - \theta \to +0$ by virtue of Lemma 2.2. Thus, the function $u(\cdot,t)$ is weakly continuous in $L_2(D)$.

Now let us prove that the function $u(\cdot,t)$ is strongly continuous to one side in $L_2(D)$. We have

$$(u(\cdot,t) - u(\cdot,\theta),\ u(\cdot,t) - u(\cdot,\theta))_{L_2(D)} = ((\mathscr{L}_{\theta,t} - I)u(\cdot,t),\ (\mathscr{L}_{\theta,t} - I)u(\cdot,t))_{L_2(D)}$$
$$= (u(\cdot,t), (\mathscr{L}^*_{\theta,t} - I)(\mathscr{L}_{\theta,t} - I)u(\cdot,t))_{L_2(D)} \to 0$$

as $\theta \to t - 0$ by virtue of Lemma 2.2.

Now we prove that the function $v(\cdot,t)$ is weakly continuous in $L_2(D)$. For $\xi \in L_2(D)$ and $\xi \leq \theta < t < T$ we have

$$(\xi, v(\cdot,t) - v(\cdot,\theta))_{L_2(D)} = (\varphi, \mathscr{L}^*_t\xi)_{X^0[t,T]} - (\varphi, \mathscr{L}^*_\theta\xi)_{X^0[t,T]} = \int_t^T (\varphi(\cdot,r), \mathscr{L}^*_{t,r}\xi - \mathscr{L}^*_{\theta,r}\xi) dr$$
$$- \int_\theta^t (\varphi(\cdot,r), \mathscr{L}^*_{\theta,r}\xi) dr = \int_t^T (\varphi(\cdot,r), \mathscr{L}^*_{t,r}(I - \mathscr{L}^*_{\theta,t})\xi)_{L_2(D)} dr - \int_\theta^t (\varphi(\cdot,r), \mathscr{L}^*_{\theta,r}\xi)_{L_2(D)} dr. \qquad (3.2)$$

Evidently, the second term on the right-hand side of (3.2) vanishes as $t - \theta \to +0$ because $\|\mathscr{L}^*_{\theta,r}\xi\|_{\mathscr{B}[\theta,T]} \leq c_0$. Besides,

$$\int_t^T (\varphi(\cdot,r), \mathscr{L}^*_{t,r}(I - \mathscr{L}^*_{\theta,r})\xi) dr = \int_t^T ((I - \mathscr{L}_{\theta,t})\mathscr{L}_{t,r}\varphi(\cdot,r), \xi)_{L_2(D)} dr.$$

For any $r \in [t,T]$ such that $\varphi(\cdot,r) \in L_2(D)$, the expression in the integral with respect to $dr$ in (3.3) vanishes in view of Lemma 2.2. Besides,

$$((I - \mathscr{L}_{\theta,t})\mathscr{L}_{t,r}\varphi(\cdot,r)\xi)_{L_2(D)} \leq c_0(1 + c_0)\|\varphi(\cdot,r)\|_{L_2(D)}\|\xi\|_{L_2(D)}.$$

Thus, we find that the quantity defined in (3.2) vanishes as $t - \theta \to +0$.

Let us prove that the function $v(\cdot,t)$ is strongly continuous to one side in $L_2(D)$. We set $\varrho = v(\cdot,t) - v(\cdot,\theta)$. For $s \leq \theta < t \leq T$, $p_1 = \mathscr{L}^*_\theta \varrho$, and $p_2 = \mathscr{L}^*_t \varrho$, we have

$$\|v(\cdot,t) - v(\cdot,\theta)\|_{L_2(D)} = (\varphi, p_2)_{X^0[t,T]} - (\varphi, p_1)_{X^0[s,T]}$$
$$= \int_t^T (\varphi(\cdot,r), p_2(\cdot,r) - p_1(\cdot,r))_{L_2(D)} dr - \int_\theta^t (\varphi(\cdot,r), p_1(\cdot,r))_{L_2(D)} dt. \quad (3.3)$$

Evidently, the second term on the right-hand side of (3.3) vanishes as $\theta \to t - 0$ because

$$\|v\|_{\mathscr{B}} \leq \text{const}, \quad \|\varrho\|_{L_2(D)} \leq \text{const} \quad (\forall \theta, t),$$
$$\|p_k(\cdot,r)\|_{L_2(D)} \leq \text{const} \quad (\forall \theta, t, r), \quad k = 1, 2. \quad (3.4)$$

Note that $p_1|_{[t,T]} = \mathscr{L}^*_t p_1(\cdot,t)$. Consequently,

$$\int_t^T (\varphi(\cdot,r), p_2(\cdot,r) - p_1(\cdot,r))_{L_2(D)} dr = \int_t^T (\varphi(\cdot,r), \mathscr{L}^*_{t,r}\varrho - \mathscr{L}^*_{\theta,r}\varrho)_{L_2(D)} dr$$
$$= \int_t^T (\varphi(\cdot,r), \mathscr{L}^*_{t,r}\varrho - \mathscr{L}^*_{t,r}\mathscr{L}^*_{\theta,r}\varrho)_{L_2(D)} dr$$
$$= \int_t^T (\mathscr{L}^*_{t,r}\varphi(\cdot,r), (I - \mathscr{L}_{\theta,r})\varrho)_{L_2(D)} dr \to 0$$

as $\theta \to t - 0$ by virtue of Lemma 2.2 and the inequality in (3.4). This completes the proof.

## 4. ESTIMATES FOR INTEGRAL FUNCTIONALS

Let $y^{x,s}(t)$ be a solution of (1.1) with the initial condition $y(s) = x$, $\tau^{x,s}(t) = \min\{\tilde{\tau}^{x,s}, T\}$, where $\tilde{\tau}^{x,s} = \inf\{t : t \geq s, y^{x,s}(t) \notin D\}$, $(x,s) \in Q$, $Q = D \times (0,T)$.

**Theorem 4.1.** *Let a function $\varphi(x,t) : \bar{D} \times [0,T] \to \mathbf{R}$ be Lebesgue measurable and $\varphi \in X^0_1$, $v = L\varphi$, and a function $V(x,s)$ be specified by (1.5), i.e.,*

$$V(x,s) = \int_s^{\tau^{x,s}} \varphi(y^{x,s}(t), t) dt.$$

*Then $V(\cdot, s) = v(\cdot, s)$ in $L_2(D)$ for all $s$ (i.e., $V(x,s) = v(x,s)$ for almost all $x$) and $V \in \mathscr{B}$.*

**Theorem 4.2.** *Let a function $\zeta(x) : \bar{D} \to \mathbf{R}$ be Lebesgue measurable and $\zeta \in L_2(D)$, $u(\cdot, s) = \mathscr{L}_{s,T}\zeta$, and a function $U(x,s)$ be specified by (1.4), i.e.,*

$$U(x,s) = \zeta(y^{x,s}(T)) \text{ Ind } \{\tilde{\tau}^{x,s} > T\}.$$

*Then $U(\cdot, s) = u(\cdot, s)$ in $L_2(D)$ for all $s$ and $U \in \mathscr{B}$.*

**Corollary 4.1.** *The functions $V(\cdot, s)$ and $U(\cdot, s)$ are weakly continuous and strongly left-continuous with respect to $s$ in $L_2(D)$. Besides,*

$$\sup_{s \in [0,T]} \|V(\cdot, s)\|_{L_2(D)} \leq c_0 \int_0^T \|\varphi(\cdot, t)\|_{L_2(D)} dt,$$

$$\sup_{s \in [0,T]} \|U(\cdot, s)\|_{L_2(D)} \leq c_0 \|\zeta\|_{L_2(D)},$$

where $c_0 = c_f(0, T)$ is defined in (1.9), i.e.,

$$c_0 = \exp \frac{1}{2} \int_0^T \sup_{x \in D} \left| \sum_{i=1}^n \frac{\partial f_i}{\partial x_i}(x, t) \right| dt.$$

*Proof of Theorem 4.1.* Let $\varrho_s \in L_2(D)$ be the distribution density of a random $n$-vector $a = a(\omega)$ such that $a \in D$ almost certainly. We assume $a = a(\omega)$ is a random vector in a probability space $(\Omega, \mathscr{F}, P)$, where $\Omega = \{\omega\} = \mathbf{R}^n$ is the Lebesgue supplement of the $\sigma$-algebra of the Borel sets in $\mathbf{R}^n$, $P$ is a probability measure on $\mathscr{F}$ such that $P(d\omega) = \varrho_s(\omega) d\omega$, in doing so we suppose $a(\omega) \equiv \omega$.

We have shown [1] that the function $\varphi(y^{a(\omega),s}(t, \omega), t) : \Omega \times [s, T] \to \mathbf{R}$ is Lebesgue measurable as a function of $\omega$ and $t$ for all $s \in [0, T]$.

By virtue of Theorem 2.1, the function $p(x, t) = \mathscr{L}_s^* \varrho_s \in \mathscr{C}^0[s, T]$ is the distribution density of the random process $y^{a,s}(t)$ defined by (1.1) and (1.2), which has a discontinuity on the boundary $\partial D$ of the domain. Hence,

$$(\varphi(\cdot, t), p(\cdot, t))_{L_2(D)} = M \left\{ \varphi(y^{a,s}(t), t) \operatorname{Ind} \{\tau^{a,s} > t\} \right\} \qquad (4.1)$$

for almost all $t$ (use Proposition 3 in [10, p. 288] and Theorem 1 in [10, p. 331]). From the definition of the operator $L$ we obtain

$$(v(\cdot, s), \varrho_s)_{L_2(D)} = \int_s^T (\varphi(\cdot, t), p(\cdot, t))_{L_2(D)} dt = \int_s^T M \left\{ \varphi(y^{a,s}(t), t) \operatorname{Ind} \{\tau^{a,s} > t\} \right\} dt. \qquad (4.2)$$

Also note that $M\{|\varphi(y^{a,s}(t), t)| \operatorname{Ind}\{\tau^{a,s} > t\}\} \in L_1(s, T)$, since we can substitute the function $\varphi$ in (4.1) by $|\varphi(x, t)|$ and derive

$$M \left\{ |\varphi(y^{a,s}(t), t)| \operatorname{Ind} \{\tau^{a,s} > t\} \right\} = (|\varphi(\cdot, t)|, p(\cdot, t))_{L_2(D)} \in L_1(s, T).$$

It follows that we can apply the Fubini theorem

$$\int_s^T M \left\{ \varphi(y^{a,s}(t), t) \operatorname{Ind} \{\tau^{a,s} > t\} \right\} dt = M \int_s^{\tau^{a,s}} \varphi(y^{a,s}(t), t) dt$$

$$= \int_\Omega P(d\omega) \int_s^{\tau^{a(\omega),s}} \varphi(y^{a(\omega),s}(t, \omega), t) dt = (V(\cdot, s), \varrho_s(\cdot))_{L_2(D)},$$

with respect to the choice of measure $P(d\omega)$ and $a(\omega) \equiv \omega$. This and (4.2) yield

$$(v(\cdot, s), \varrho_s)_{L_2(D)} = (V(\cdot, s), \varrho_s)_{L_2(D)} \qquad (4.3)$$

for all $s \in [0, T]$, $\varrho_s \in L_2(D)$ [such that $\varrho_s(x) \geq 0$ and $\int_D \varrho_s(x) dx = 1$]. Then we have proved (4.3) for an arbitrary $\varrho_s \in L_2(D)$ using standard reasoning. Consequently, $V(\cdot, s) = v(\cdot, s)$ in $L_2(D)$ for all $s \in [0, T]$. Theorem 4.1 is proved.

*Proof of Theorem 4.2.* Let $(\Omega, \mathscr{F}, P)$, $a = a(\omega)$, $\varrho_s \in L_2(D)$ be the same as in the proof of Theorem 4.1, $\Omega = \{\omega\} = D$.

From the proof of Theorem 1.1 [1], the function $\zeta(y^{a(\omega),s}(T, \omega)) : \Omega \to \mathbf{R}$ is Lebesgue measurable, as is the function Ind $\{\tilde{\tau}^{a(\omega),s} > T\} : \Omega \to \mathbf{R}$. We have

$$(u(\cdot,s), \varrho_s)_{L_2(D)} = (\zeta, \mathscr{L}_{s,T}\varrho_s)_{L_2(D)} = M\{\zeta(y^{a,s}(T))\operatorname{Ind}\{\tilde{\tau}^{a,s} > T\}\} = (V(\cdot,s), \varrho_s)_{L_2(D)}.$$

From this we obtain the desired result in the same manner as for the proof of Theorem 4.1. This completes the proof of Theorem 4.2.

Finally, we note that results similar to those in this paper were announced in [11], where we considered general degenerate Itô equations. We also mention [12], where we derived estimates like those in [1] for Itô processes, and [13], where nondegenerate Itô processes in a bounded domain are studied.

The author is grateful to S. Yu. Pilyugin for his useful consultations and the idea behind the proof of Proposition 2.1.